\documentclass[a4paper,8pt]{article}
\usepackage{amsmath}
\usepackage{amsthm,amsfonts,amssymb}
\usepackage[dvips]{graphicx}
\usepackage[english]{babel}
\usepackage{color}
\usepackage{graphics}
\usepackage{graphpap}

\providecommand{\U}[1]{\protect\rule{.1in}{.1in}}

\providecommand{\U}[1]{\protect\rule{.1in}{.1in}}
\newtheorem{theorem}{Theorem}[section]

\newtheorem{corollary}[theorem]{Corollary}
\newtheorem{definition}[theorem]{Definition}
\newtheorem{lemma}[theorem]{Lemma}
\newtheorem{proposition}[theorem]{Proposition}

\makeatletter
\def\blfootnote{\xdef\@thefnmark{}\@footnotetext}
\makeatother

\begin{document}

\title{Transient analysis of a $M/M/\infty$ queue with discouragement and for the related embedded chain}
\author{Andrea Monsellato}
\maketitle

\begin{abstract}
\noindent
Consider the following birth and death process with the following infinitesimal transition probabilities $\lambda_k=\frac{\lambda}{1+k}$ and
$\mu_k=\mu k$ with $\lambda,\mu>0$.
\noindent
This process has known as a \emph{discouragement queue} \cite{Natvigtransientdisc}.\\
\noindent
Although from the theoretical point of view the problem of the determination of the transition functions has been solved \cite{Karlindbdp}, the explicit form of them is not present in literature.\\
\noindent
We have solved this problem assuming that the solution is representable by a Taylor series, under the initial condition that the process starts to state zero.\\
\noindent
We discuss also the same problem for the embedded chain and using direct computation we obtain a recursive formula for the transient distribution.\\
\end{abstract}


\section{Introduction}

\noindent
The problem to find an explicit formula of the transition probabilities for this particular queue is proposed by Natvig \cite{Natvigtransientdisc}.\\
\noindent
Parthasarathy et al. \cite{Partha} give an explicit solution, based on continued fractions approach, for the cases $\lambda_k=\frac{\lambda}{1+k}$ and $\mu_k=k$ and $\lambda_k=\lambda$ and $\mu_k=\mu k$, with $\lambda,\mu>0$.\\
\noindent
Despite the powerful of continued fractions approach, for our case this method seems to be fail.\\
\noindent
We propose a more simple approach based on the Taylor expansion, and recognize a simple iterative formula for the coefficients of the series.
Our approach permits a fast numerically calculation of the coefficients.\\
\noindent
Unexpectedly the bound we used for the coefficients, to verify the total convergence of the Taylor series, is related with the Bessel Numbers \cite{bessel}.\\

\section{Transient solution in continuous time}

\noindent
We remember that that the infinitesimal transition probabilities are given by

\begin{equation*}
P(X(t+h)=n+m|X(t)=n)=\left\{ \begin{matrix} \frac{\lambda}{1+k}h+o(h), &\mbox{if}\quad m=1\\ o(h),  &\mbox{if}\quad |m|>1 \\\mu kh+o(h), &\mbox{if} \quad m=-1\end{matrix}\right.
\end{equation*}

\noindent
Let $q(k,t)=P(X(t)=k)$ then

\begin{align}\label{qkbdp}
\left\{
\begin{array}{ll}
q'(0,t)=-\lambda q(0,t)+\mu q(1,t)\\
q'(k,t)=\frac{\lambda}{k}q(k-1,t)-\left(\frac{\lambda}{1+k}+k\mu \right)q(k,t)+(k+1)\mu q(k+1,t),\quad \forall k\geq 1
\end{array}
\right.
\end{align}

\noindent
Let $\tau=\lambda t$, $p(k,\tau)=q(k,t(\tau))$, $\alpha=\sqrt{\frac{\mu}{\lambda}}$ and
$\dot{p}(k,\tau)=\frac{1}{\lambda}\dot{q}(k,t)$, then (\ref{qkbdp}) is equivalent to

\begin{align}\label{pkbdp}
\begin{small}
\left\{
\begin{array}{ll}
\dot{p}(0,\tau)=-p(0,\tau)+\alpha^2 p(1,\tau)\\
\dot{p}(k,\tau)=\frac{1}{k}p(k-1,\tau)-\left(\frac{1}{1+k}+k\alpha^2 \right)p(k,\tau)+(k+1)\alpha^2 p(k+1,\tau),\quad \forall k\geq 1
\end{array}
\right.
\end{small}
\end{align}

\bigskip
\noindent
We consider (\ref{pkbdp}) under the initial conditions

\begin{align*}
p(0,0)=1, \quad p(k,0)=0 \quad \forall k>0
\end{align*}

\bigskip
\noindent
Let us observe that the functions $p(k,\cdot)$ have derivatives of any orders. This observations suggest a Taylor series expansion for
$p(k,\tau)$, $k\geq 0$, and in particular for $p(0,\tau)$.\\

\noindent
For sake of simplicity we put $p(k,\tau)=\frac{w(k,\tau)}{\alpha^k k!}$, then the $p(k,\tau)$, $\tau \geq 0$, satisfy (\ref{pkbdp})
if and only if the $w(k,\tau)$, $k\geq 0$, satisfy

\begin{align}\label{wkbdp}
\begin{small}
\left\{
\begin{array}{ll}
\dot{w}(0,\tau)=-w(0,\tau)+\alpha w(1,\tau)\\
\dot{w}(k,\tau)=\alpha w(k-1,\tau)-\left(\frac{1}{1+k}+k\alpha^2\right)w(k,\tau)+\alpha w(k+1,\tau),\quad \forall k\geq 1
\end{array}
\right.
\end{small}
\end{align}

\bigskip
\noindent
with the initial conditions $w(0,0)=1$, $w(k,0)=0$.\\

\noindent
Let $b_k=\frac{1}{k+1}+k\alpha^2$, $k\geq 0$, we observe that if

\begin{align}\label{taylorsolutionbdp}
&w(0,\tau)=\sum_{i=0}^{+\infty}\frac{\tau^i}{i!}r_i,\quad r_0=1
\end{align}

\bigskip
\noindent
then

\begin{align*}
&w(1,\tau)=\frac{1}{\alpha}\sum_{i=0}^{+\infty}\frac{\tau^i}{i!}[r_{i+1}+b_0r_i];\\&
\phantom{x}\\&
w(2,\tau)=\frac{1}{\alpha^2}\sum_{i=0}^{+\infty}\frac{\tau^i}{i!}[r_{i+2}+b_0r_{i+1}+b_1(r_{i+1}+b_0r_i)-\alpha^2 r_i];\\&
\phantom{x}\\&
w(3,\tau)=\frac{1}{\alpha^3}\sum_{i=0}^{+\infty}\frac{\tau^i}{i!}\big\{(r_{i+3}+b_0r_{i+2})+b_1(r_{i+2}+b_0r_{i+1})-\alpha^2 r_{i+1}+\\&
+b_2[r_{i+2}+b_0r_{i+1}+b_1(r_{i+1}+b_0r_i)-\alpha^2 r_i]-\alpha^2(r_{i+1}+b_0r_i)\big\};\\&
...
\end{align*}

\bigskip
\noindent
So that if we put $S^{(0)}_i=r_i$, $S^{(-1)}_i=0$ for all $i\geq 0$ and $S^{(n)}_i=0$ for all $n>i$, by recursion we have

\begin{align}\label{skcoefftaylor}
\left\{
\begin{array}{ll}
S^{(1)}_i=S^{(0)}_{i+1}+b_0 S^{(0)}_i\\
S^{(k+1)}_i=S^{(k)}_{i+1}+b_k S^{(k)}_i-\alpha^2 S^{(k-1)}_i,\quad \forall i\geq 0, \quad 0\leq k\leq i+1
\end{array}
\right.
\end{align}

\bigskip
\noindent
so that

\begin{align}\label{wbdseriessol}
w(k,\tau)=\frac{1}{\alpha^k}\sum_{i=0}^{+\infty}\frac{\tau^i}{i!}S^{(k)}_{i},\quad \forall k\geq 0
\end{align}

\bigskip
\noindent
and obviously

\begin{align}\label{bdseriessol}
p(k,\tau)=\frac{1}{\alpha^{2k} k!}\sum_{i=0}^{+\infty}\frac{\tau^i}{i!}S^{(k)}_{i},\quad \forall k\geq 0
\end{align}

\bigskip

\begin{proposition}
Let $S^{(0)}_{0}=1$ and $S^{(-1)}_{i}=1$, for all $i\geq 0$, then (\ref{skcoefftaylor}) is equivalent to

\begin{align}\label{sncoeffsum}
\left\{
\begin{array}{ll}
 S^{(k)}_k=1, \quad \forall k \geq 0\\
S^{(k)}_{h+k+1}=\alpha^2 S^{(k-1)}_{k+h}-b_k S^{(k)}_{k+h}+S^{(k+1)}_{k+h}, \quad \forall h,k\geq 0
\end{array}
\right.
\end{align}

\bigskip
\noindent
furthermore

\begin{align}\label{sncoeffsum2}
S^{(k)}_{k+h+1}=\sum_{i=0}^{k}\alpha^{2(k-i)}\left[ -b_i S^{(i)}_{i+h}+S^{(i+1)}_{i+h}\right], \quad h,k\geq 0
\end{align}

\noindent
\textbf{Proof}\\

\end{proposition}

\noindent
(\ref{sncoeffsum}) is obvious. To prove (\ref{sncoeffsum2}) we proceed by induction.\\
\noindent
Thesis is true for $n=0$ because

\begin{align*}
S^{(0)}_{h+1}=\alpha^2 S^{(-1)}_{h}-b_0 S^{(0)}_{h}+S^{(1)}_{h}=\sum_{i=0}^{0}\alpha^{-2i}\left[ -b_i S^{(i)}_{i+h}+S^{(i+1)}_{i+h}\right]
\end{align*}

\bigskip
\noindent
then supposing it is true for $k>0$ we have

\begin{align*}
S^{(k+1)}_{(k+1)+h+1}&=\alpha^2 S^{(k+1-1)}_{k+1+h}-b_{k+1} S^{(k+1)}_{k+1+h}+S^{(k+1+1)}_{k+1+h}=\\&
\quad\\&
=\alpha^2 S^{(k)}_{k+1+h}-b_{k+1} S^{(k+1)}_{k+1+h}+S^{(k+2)}_{k+1+h}=\\&
\mbox{}\\&
=\alpha^2\sum_{i=0}^{k}\alpha^{2(k-i)}\left[ -b_i S^{(i)}_{i+h}+S^{(i+1)}_{i+h}\right]+\left[-b_{k+1} S^{(k+1)}_{k+1+h}+S^{(k+2)}_{k+1+h}\right]=\\&
\quad\\&
=\sum_{i=0}^{k}\alpha^{2(k+1-i)}\left[ -b_i S^{(i)}_{i+h}+S^{(i+1)}_{i+h}\right]+\alpha^{2(k+1-(k+1))}\left[-b_{k+1} S^{(k+1)}_{k+1+h}+S^{((k+1)+1)}_{k+1+h}\right]=\\&
\quad\\&
=\sum_{i=0}^{k+1}\alpha^{2(k+1-i)}\left[ -b_i S^{(i)}_{i+h}+S^{(i+1)}_{i+h}\right]
\end{align*}

\begin{flushright}
$\Box$
\end{flushright}

\noindent
To establish the convergence of the Taylor series (\ref{bdseriessol}) we observe that if we put

\begin{align}\label{lkbdpsol}
& L^{(k)}_i\doteq \frac{(-1)^{i-k}}{\alpha^{2k}} S^{(k)}_i
\end{align}

\noindent
from (\ref{skcoefftaylor}) we have

\begin{align*}
&L^{(k)}_{i+1}=L^{(k-1)}_{i}+b_k L^{(k)}_i+\alpha^2 L^{(k+1)}_i, \quad \forall i\geq 0, \quad 0\leq k\leq i+1
\end{align*}

\noindent
and

\begin{align*}
&L^{(-1)}_{0}=0,\quad \forall i\geq 0;\quad L^{(k)}_{i}=0,\quad \forall k>i;\quad L^{(i)}_{i}=1,\quad \forall i\geq 0
\end{align*}

\bigskip
\noindent
It is obvious that such $L^{(k)}_{i}$ are non negative for all $i,k=0,1,2,...$ and moreover if we put $\gamma=\max(1,\alpha^2)$ we have

\begin{align}\label{lkipositive}
&L^{(k)}_{i}\leq M^{(k)}_{i} \gamma^{2(i-k)},\quad \forall k,i=0,1,2,...
\end{align}

\bigskip
\noindent
where the $M^{(k)}_{i}$ are such that

\begin{align}\label{besseln}
\left\{
\begin{array}{ll}
M^{(k)}_{i+1}=M^{(k-1)}_{i}+(1+k) M^{(k)}_i+M^{(k+1)}_i, \quad \forall i\geq 0, \quad 0\leq k\leq i+1\\
M^{(-1)}_{0}=0,\quad \forall i\geq 0;\quad M^{(k)}_{i}=0,\quad \forall k>i;\quad M^{(i)}_{i}=1,\quad \forall i\geq 0
\end{array}
\right.
\end{align}

\newpage
\noindent
Computing the $M^{(k)}_{i}$, we obtain that

\begin{align*}
M^{(0)}_{i}=B^*_i, \quad \forall i\geq 0
\end{align*}

\bigskip
\noindent
where the $B^*_i$ are the so called Bessel numbers. This last assertion derives from the comparison of the generating function related to the solution of the equations (\ref{besseln}), see \cite{aigner} for more details (paragraph 7.3), and the generating function of the Bessel Number proposed in \cite{bessel}.\\

\noindent
In \cite{bessel} it is proved that Bessel numbers have the asymptotic form

\begin{align*}
B^*_i\sim \frac{1}{\sqrt{2\pi i}}\frac{w^{i+3}}{(w!)^2}, \quad \forall i\geq 0
\end{align*}

\bigskip
\noindent
where $w\sim \frac{i}{2\ln( i)}$ is the positive root of following equation:

\begin{align*}
i+2=2w\ln( w)
\end{align*}

\bigskip
\noindent
Then, apart from sub-exponential factors, Bessel numbers grown like

\begin{align}\label{bellasym}
B^*_i\approx \left(\frac{i}{2e\ln(i)}\right)^i
\end{align}

\bigskip
\noindent
From (\ref{bellasym}) it follows that the power series

\begin{align}\label{besselseriesupp}
\sum_{i=0}^{+\infty}(-1)^i\frac{\tau^i}{i!}M^{(0)}_{i}
\end{align}

\bigskip
\noindent
has an infinite convergence radius.\\

\noindent
Therefore taking into account (\ref{lkipositive}) by the convergence of (\ref{besselseriesupp}), we conclude that the series at the second member of (\ref{wbdseriessol}) is a power series with infinite convergence radius.\\

\noindent
In conclusion we have established the following theorem.

\begin{theorem}
The Cauchy problem (\ref{taylorsolutionbdp}) has a unique solution $p(k,\tau)$, $k=0,1,2,...$, where $p(k,\tau)$, $\tau \geq 0$, are defined in (\ref{bdseriessol}).
\end{theorem}

\newpage
\section{Embedded chain}

\noindent
We start with the general state dependent case and we recover the solution of our problem as corollary of the general case.\\
\noindent
We consider now the following birth and death process with infinitesimal transition probabilities

\begin{equation*}
P(X(t+h)=k+m|X(t)=k)=\left\{ \begin{matrix} \lambda_k h+o(h), &\mbox{if}\quad m=1\\ o(h),  &\mbox{if}\quad |m|>1 \\\mu_k h+o(h), &\mbox{if} \quad m=-1\end{matrix}\right.
\end{equation*}

\noindent
Let $p(k,t)=P(X(t)=k)$ then

\begin{align}\label{qkbdp}
\left\{
\begin{array}{ll}
p'(0,t)=-\lambda p(0,t)+\mu p(1,t)\\
p'(k,t)=\lambda_{k-1}p(k-1,t)-(\lambda_k+\mu_k)p(k,t)+\mu_{k+1} p(k+1,t),\quad \forall k\geq 1
\end{array}
\right.
\end{align}

\bigskip

\begin{definition}\label{def:regularjumpprocess}
A stochastic process $\{X(t)\}_{t\geq 0}$ taking its values in the countable state space $E$ is a called \emph{jump process}
if for the almost $\omega\in\Omega$ and all $t\geq 0$, there exists $\epsilon(t,\omega)>0$ such that

\begin{equation*}
X(t+s,\omega)=X(t,\omega),\quad \forall s\in[t,t+\epsilon(t,\omega))
\end{equation*}

\bigskip
\noindent
It is called \emph{regular jump process} if in addiction, for almost all $\omega\in\Omega$, the set $A(\omega)$ of
discontinuities of the function $t\rightarrow X(t,\omega)$ is $\sigma-$discrete, that is, for all $c\geq 0$

\begin{equation*}
|A(\omega)\cap[0,c]|<+\infty
\end{equation*}

\bigskip
\noindent
where the notation $|B|$ is the cardinality of set $|B|$. A \emph{regular jump homogeneus Markov chain} is by definition
a continuous time \textbf{HMC} that is also regular jump process.
\end{definition}

\noindent
Let $\{\tau_n\}$ be a non decreasing sequence of transition times of the regular jump process $\{X(t)\}_{t\geq 0}$
where $\tau_0=0$ and $\tau_n=\infty$ if there are strictly fewer than $n$ transitions in $(0,\infty)$.\\

\noindent
The process $\{X_n\}_{n\geq 0}$ with value in $E_{\Delta}=E\cup\Delta$, where $\Delta$ is an arbitrary element not in $E$,
is defined by $X_n=X(\tau_n)$ with the convention $X(\infty)=\Delta$, and it is called \emph{embedded process} of the jump process.\\

\noindent
The associated embedded process, see \cite{queueingtheorygross}, has transition probabilities of the form

\begin{align*}
\left\{
\begin{array}{ll}
p_{i,j}=\frac{\lambda_i}{\lambda_i+\mu_i}, \quad j=i+1, \quad i\geq 1\\
p_{i,j}=\frac{\mu_i}{\lambda_i+\mu_i}, \quad j=i-1, \quad i\geq 1\\
p_{i,j}=1, \quad j=1, \quad i=0\\
p_{i,j}=0 \quad otherwise
\end{array}
\right.
\end{align*}

\noindent
Now consider the associated embedded processes $\{X_n\}_{n\geq 0}$ of (\ref{qkbdp}).\\
\noindent
Let $P(X_n=k)=p_{n,k}$ then

\begin{align}\label{equ:bddiscretamod}
\left\{
\begin{array}{ll}
p_{n+1,k}=\frac{\lambda_{k-1}}{\lambda_{k-1}+\mu_{k-1}}p_{n,k-1}+\frac{\mu_{k+1}}{\lambda_{k+1}+\mu_{k+1}}p_{n,k+1},\quad \forall n\geq 0,
\quad 1\leq k\leq n+1\\
p_{n,k}=0,\quad \forall k>n\\
p_{n,-1}=0,\quad \forall n\geq 0\\
p_{0,0}=1
\end{array}
\right.
\end{align}

\noindent
The next two following lemma are elementary and they will be explained without proofs.

\begin{lemma}\label{pnkrecbdpdis}
\noindent
Let $(p_{n,k})_{n,k\geq 0}$ be the matrix defined in (\ref{equ:bddiscretamod}) then

\begin{align*}
&(i)\quad p_{n,n}=\prod_{i=0}^{n}\left(\frac{\lambda_i}{\lambda_i+\mu_i}\right), \quad \forall n> 0\\&
\phantom{}\\&
(ii)\quad \sum_{k=0}^{+\infty}p_{n,k}=1, \quad \forall n\geq 0
\end{align*}

\end{lemma}

\noindent
From (\ref{equ:bddiscretamod}) it follows obviously the following lemma.\\

\begin{lemma}\label{lem:bdparity}
\noindent
In the hypothesis of lemma (\ref{pnkrecbdpdis}), if $n+k$ is odd  then $p_{n,k}=0$.
\end{lemma}

\noindent
By virtue of Lemma (\ref{pnkrecbdpdis}) and Lemma (\ref{lem:bdparity}) we can find the transient distribution.

\begin{proposition}\label{prop:bdtransientsolution}

\noindent
Let $(p_{n,k})_{n,k\geq 0}$ as in (\ref{equ:bddiscretamod}) then

\begin{align}\label{equ:transientbdsolution}
\left\{
\begin{array}{ll}
p_{n,k}=\prod_{i=0}^{k}\left(\frac{\lambda_i}{\lambda_i+\mu_i}\right)T_{\frac{n-k}{2},k},\quad n\geq 0,\quad 0\leq k\leq n,\quad n+k\quad even\\
p_{n,k}=0,\quad otherwise
\end{array}
\right.
\end{align}

\noindent
where

\begin{align*}
&T_{h,k}=\sum_{l=0}^{k}\frac{\lambda_l}{\lambda_l+\mu_l}\left(1-\frac{\lambda_{l+1}}{\lambda_{l+1}+\mu_{l+1}}\right)T_{h-1,l+1},\quad \forall k\geq 0,h\geq 1;\quad T_{0,k}=1,\quad \forall k\geq 0
\end{align*}

\bigskip
\noindent
\textbf{Proof}\\

\end{proposition}

\noindent
We know from the lemma (\ref{lem:bdparity}) that if $n+k$ is odd then $p_{n,k}=0$, furthermore it holds that

\begin{equation}\label{dkparityproduct}
d_k:=p_{k,k}=\prod_{i=0}^{k}\left(\frac{\lambda_i}{\lambda_i+\mu_i}\right)
\end{equation}

\noindent
Let $n=k+2h$, we want prove that

\begin{equation}\label{equ:bdparitysol}
p_{k+2h,k}=d_k T_{h,k}, \quad \forall k\geq 0,h\geq 0
\end{equation}

\noindent
Fixed $k$, we proceed by induction on $h$. The (\ref{equ:bdparitysol}) is true for $h=0$ being

\begin{equation*}
p_{k,k}=d_k=d_k T_{0,k}
\end{equation*}

\noindent
supposing it is true for $h$ we prove it for $h+1$. Let $\alpha_k:=\frac{\lambda_k}{\lambda_k+\mu_k}$ from (\ref{equ:bddiscretamod}) we have

\begin{align*}
p_{k+2(h+1),k}&=p_{k+2h+1,k-1}\alpha_{k-1}+p_{k+2h+1,k+1}\alpha_{k+1}=\\&\\&
=p_{k+2h+1,k-1}\alpha_{k-1}+[1-\alpha_{k+1}]p_{k+2h+1,k+1}
\end{align*}

\noindent
and by the inductive hypothesis (\ref{equ:bdparitysol}), it follows that

\begin{align*}
&p_{k+2(h+1),k}=p_{k+2h+1,k-1}\alpha_{k-1}+[1-\alpha_{k+1}]d_{k+1}T_{h,k+1}
\end{align*}

\noindent
Taking into account (\ref{dkparityproduct}) we have

\begin{align}\label{equ:bdpassoinduttivo}
&p_{k+2(h+1),k}=p_{k+2h+1,k-1}\alpha_{k-1}+(d_{k+1}-d_{k+2})T_{h,k+1}
\end{align}

\noindent
On other hand for every fixed $h$, setting

\begin{equation}\label{equ:bhktransbdsolution}
b^h_k=p_{k+2h+1,k-1}\alpha_{k-1}
\end{equation}

\noindent
we have (to show)

\begin{align}\label{equ:bhktransbdsolutionfinal}
\left\{
\begin{array}{ll}
b^h_0=0\\
b^h_k=d_k\sum_{i=0}^{k-1}\frac{d_{i+1}-d_{i+2}}{d_i}T_{h,i+1},\quad \forall k\geq 1
\end{array}
\right.
\end{align}

\noindent
In fact proceeding by induction also in this case, we have that (\ref{equ:bhktransbdsolutionfinal}) is true for $k=0$, being  $p_{2h+1,-1}=0$.
Furthermore it is true for $k=1$ because for (\ref{equ:bhktransbdsolution}) it holds that

\begin{equation*}
b^h_1=p_{2+2h,0}
\end{equation*}

\noindent
and from (\ref{equ:bddiscretamod}) we have

\begin{equation*}
b^h_1=p_{1+2h,-1}+(1-\alpha_1)p_{1+2h,1}=(1-\alpha_1)p_{1+2h,1}
\end{equation*}

\noindent
thus by inductive hypothesis on $h$ it follows that

\begin{equation*}
b^h_1=(1-\alpha_1)d_1T_{h,1}=(d_1-d_2)T_{h,1}=d_1\frac{d_1-d_2}{d_0}T_{h,1}
\end{equation*}

\noindent
being $d_0=d_1=1$.\\

\noindent
Now supposing that (\ref{equ:bhktransbdsolutionfinal}) is true for $k\geq 1$, we prove it for $k+1$:

\begin{align*}
b_{k+1}^h&=p_{k+1+2h+1,k}\alpha_{k}=\\&\\&
=\alpha_{k}\left[p_{k+1+2h,k}\alpha_{k}+(1-\alpha_{k+1})p_{k+1+2h,k+1}\right]
\end{align*}

\noindent
for the (\ref{equ:bhktransbdsolution}) and by inductive hypothesis on  $h$ (\ref{equ:bdparitysol}) we obtain

\begin{equation*}
b_{k+1}^h=\alpha_{k}[b^h_k+(1-\alpha_{k+1})d_{k+1}T_{h,k+1}]
\end{equation*}

\noindent
furthermore by inductive hypothesis on $k$ (\ref{equ:bhktransbdsolutionfinal}) we have

\noindent
\begin{align*}
b_{k+1}^h&=\alpha_{k}\left[d_k\sum_{i=0}^{k-1}\frac{d_{i+1}-d_{i+2}}{d_i}T_{h,i+1}+
(1-\alpha_{k+1})d_{k+1}T_{h,k+1}\right]=\\&\\&
=\alpha_{k}\left[d_{k}\sum_{i=0}^{k-1}\frac{d_{i+1}-d_{i+2}}{d_i}T_{h,i+1}+
(d_{k+1}-d_{k+2})T_{h,k+1}\right]=\\&\\&
=\alpha_{k}d_{k}\sum_{i=0}^{k}\frac{d_{i+1}-d_{i+2}}{d_i}T_{h,i+1}=\\&\\&
=d_{k+1}\sum_{i=0}^{k}\frac{d_{i+1}-d_{i+2}}{d_i}T_{h,i+1}
\end{align*}

\noindent
Now for $h$ fixed the (\ref{equ:bhktransbdsolution}) is true, substituting (\ref{equ:bhktransbdsolution}) in
(\ref{equ:bdpassoinduttivo}) finally we obtain

\begin{align*}
p_{k+2(h+1),k}&=d_{k+1}\sum_{i=0}^{k}\frac{d_{i+1}-d_{i+2}}{d_i}T_{h,i+1}+
(d_{k+1}-d_{k+2})T_{h,k+1}=\\&\phantom{}\\&
=d_{k+1}\left[\sum_{i=0}^{k}\frac{d_{i+1}-d_{i+2}}{d_i}T_{h,i+1}+
\frac{d_{k+1}-d_{k+2}}{d_{k+1}}T_{h,k+1}\right]=\\&\phantom{}\\&
=d_{k+1}\sum_{i=0}^{k+1}\frac{d_{i+1}-d_{i+2}}{d_i}T^{(h)}_{i+1}=d_{k+1}T_{h+1,k}
\end{align*}

\noindent
then the thesis follows.

\begin{flushright}
$\Box$
\end{flushright}

\noindent
As an immediate corollary of the proposition \ref{prop:bdtransientsolution} we have the following:

\begin{corollary}

\noindent
Let $(p_{n,k})_{n,k\geq 0}$ as in (\ref{equ:bddiscretamod}) with $\lambda_k=\frac{\lambda}{1+k}$ and
$\mu_k=\mu k 1_{\{k\geq 1\}}$ with $\lambda,\mu>0$, then

\begin{align}\label{equ:transientbdsolution}
\left\{
\begin{array}{ll}
p_{n,k}=d_kT^{\left(\frac{n-k}{2}\right)}_k,\quad n\geq 0,\quad 0\leq k\leq n,\quad n+k\quad even\\
p_{n,k}=0,\quad otherwise
\end{array}
\right.
\end{align}

\noindent
where

\begin{align*}
&d_k=\prod_{i=0}^{k}\frac{1}{1+i(i-1)\alpha^2}\\&
T^{(h)}_k=\sum_{i=0}^{k}\frac{d_{i+1}-d_{i+2}}{d_i}T^{(h-1)}_{i+1},\quad \forall k\geq 0,h\geq 1;\quad T^{(0)}_k=1,\quad \forall k\geq 0
\end{align*}

\end{corollary}


\begin{thebibliography}{10}
\bibitem[1]{aigner}Martin A., A course in enumeration, First Edition, Springer, 2010
\bibitem[2]{bessel}Flajolet, P., Schott, R., Nonoverlapping partitions, continued fractions, {B}essel
              functions and a divergent series, European J. Combin. 11, (1990), 421--432
\bibitem[3]{queueingtheorygross} Gross, D., Harris, C.M., Fundamentals of queueing theory, Second Edition,
              John Wiley \& Sons Inc., New York, 1985
\bibitem[4]{Karlindbdp} Karlin, S., McGregor, J. L., The differential equations of birth-and-death processes, and
              the {S}tieltjes moment problem, Trans. Amer. Math. Soc. 85, (1957), 489--546
\bibitem[5]{Natvigtransientdisc} Natvig, B., On the transient state probabilities for a queueing model
              where potential customers are discouraged by queue length, J. Appl. Probability 11, 1974, 345--354
\bibitem[6]{Partha} Parthasarathy, P.R., Selvaraju, N., Transient analysis of a queue where potential customers are discouraged by queue length,
                    Mathematical Problems in Engineering 7, 5 (2001), 433-454
\end{thebibliography}

\end{document}